\documentclass{article}
\title{Formation of Multi-Point Singularities of Self-Similar Type for Burgers Equation}
\date{}
\author{Yiya Qiu\footnote{School of Mathematical Sciences, University of Science and Technology of China, Hefei, Anhui, 230026, PR China, Email:qq171579@mail.ustc.edu.cn}, Lifeng Zhao\footnote{School of Mathematical Sciences, University of Science and Technology of China, Hefei, Anhui, 230026, PR China, Email:zhaolf@ustc.edu.cn}}

\makeatletter
\@addtoreset{equation}{section}
\makeatother

\usepackage{xcolor}

\usepackage{amsmath}
\usepackage{amssymb}
\usepackage{amsthm}
\usepackage{mathtools}
\usepackage[a4paper,left=3cm,right=3cm,top=3cm,bottom=3cm]{geometry}
\usepackage{subeqnarray}
\usepackage{cases}
\usepackage{setspace}

\newtheorem{thm}{Theorem}[section]
\newtheorem{lem}[thm]{Lemma}
\newtheorem{cor}[thm]{Corollary}
\newtheorem{prop}[thm]{Proposition}

\theoremstyle{remark}

\newcommand\al{\alpha} 
\newcommand\ve{\varepsilon}  
\newcommand\de{\delta}   
\newcommand\ka{\kappa} 

\newcommand\les{\lesssim} 
\newcommand\gtr{\gtrsim}  
\newcommand\f{\frac}  
\newcommand\p{\partial}  

\newcommand\td{\tilde} 
\begin{document}
\maketitle
\begin{abstract}
In this paper, we constuct the multi-point blowup solutions of self-similar type for the inviscid Burgers equation. The shape and blowup dynamics are precisely described. Moreover, the solutions we construct are stable under small perturbations on initial data restricted in a compact set.
\end{abstract}
\bibliographystyle{plain}
\section{Introduction}

We consider the inviscid Burgers equation on $\mathbb{R}$,
\begin{equation}\label{bur}
\begin{cases}
u_t+uu_x=0,\\
u(T-\de,x)=u_0(x),
\end{cases}
\end{equation}
where $u: [T-\de,T)\times\mathbb{R} \to\mathbb{R}$ is a real-valued function, $T>0$ is the time when we expect $u$ blow up, and $0<\de\ll1$ is under determined. 
First introduced by Bateman and Burgers in their independent works \cite{BATEMAN1915} and \cite{Burgers1929}, the equation (\ref{bur}) was studied as the simplest quasilinear model of conservation laws, nonlinear wave equations and compressible fluid dynamics, describing the generic blowup dynamics of hyperbolic systems, see for \cite{Dafermos2005} and \cite{Serre1999}.

 It is universally acknowledged that Burgers equation admits singular solutions starting from small smooth data, as long as its slope has a negative point. This is illustrated by the method of characteristics, namely the solution becomes singular in finite time once there exist two characteristic lines getting across, see \cite{Eggers2015}. At this time, the gradient of solution grows to infinity instead of the solution itself. In particular, this type of blow up is called ``geometric type''.
The study of geometric type of singularities in generic hypobolic systems is flourishing in past decades. For example, by using characteristic curves or surfaces, \cite{Alinhac1999} and \cite {Alinhac1999a} established shock formation theory on quasilinear wave equations in dimension two, \cite{Speck2014} on quasilinear wave equations on dimension three in relativistic scheme, as well as \cite{Chris2007} on compressible Euler equation in dimension three. 

Among various shock waves, Self-similar phenomena  are observed in massive physical researches and activities, which reveal the processing that scales reproducing themselves smaller and smaller. In particular, smooth self-similar profiles are always the asymptotic attractors of all smooth and non-degenerate shocks. The self-similar profile of Burgers equation was studied in \cite{Eggers2008} as an important example. And in \cite{Collot2018}, the authors summarized that 
 (\ref{bur}) admits a family of self-similar solutions
 \begin{equation} \label{selfu}
 u=(T-t)^{\f{1}{2i}}\Psi_i\left(\f{x}{(T-t)^{1+\f{1}{2i}}}\right),\ \ \ \ i\in\mathbb{N},
\end{equation}
 where $\Psi_i(X)$ is a family of odd, decreasing, analytic profiles in  $X=\f{x}{(T-t)^{1+\f{1}{2i}}}$. In particular, each $\Psi_i(X)$ satisfies the self-similar equation 
\begin{equation}\label{sebur}
-\f{1}{2i}\Psi_i+\f{2i+1}{2i}X\p_X\Psi_i+\Psi_i\p_X\Psi_i=0.
\end{equation}
Therefore, for every $i$, (\ref{selfu}) is a singular self-similar solution to (\ref{bur}), blowing up at the single point $x=0$ and the finite time $t=T$. There are also several recent works considering singular solutions of self-similar type, such as \cite{Buckmaster2019a} on 2d compressible Euler equation with azimuthal symmetry and \cite{Buckmaster2019} on 3d compressible Euler equation without symmetry assumptions.

In this article, our goal is to find a solution blowing up at finite many points at the same time of the following form:
\begin{equation}\label{deco1}
u=\sum_{l=1}^{L}u_l+\td{\ve},
\end{equation}
where 
\begin{equation}
u_l=(T-t)^{\f{1}{2i_l}}\Psi_{i_l}\left(\f{x-y_l(t)}{(T-t)^{\f{2i_l+1}{2i_l}}}\right)\end{equation}
with $y_1(t)=0$. $\td{\ve}$ is a small perturbation whose initial data $\td{\ve}_0$ is supported in 
\begin{equation}\label{sptx}
\text{supp}\ \td{\ve}_0\subseteq \left\{x: \f{4}{3}\de^{\f{1}{2i_1}}\leq x \leq \f{5}{3}\de^{\f{1}{2i_1}} \right\}.
\end{equation}
$u_l$ can be regarded as $u_1$ shifting a distance of $y_l(t)$, and we assume
\begin{equation}\label{yl0}
y_{l,0}=y_l(T-\de)=3(l-1)\de^{\f{1}{2i_1}}.
\end{equation}

Now we state the main result roughly and the precise statement is in Theorem 4.1.

\begin{thm}
There exists $\de>0$ and an open set of a suitable topology such that if $\td{\ve}_0$ is in this open set, then the equation \emph{(\ref{bur})} admits a solution of the form \emph{(\ref{deco1})}.
\end{thm}

We shall explain our main ideas. First, since our task is to construct singular solutions with finite many self-similar profiles, the interaction between different profiles is of vital importance. So we introduce $y_{l}(t),$ measuring the distance between the two centers of the first and the $l^{th}$ profiles, and figure it out by using characteristic lines. Next, we assume the initial data of disturbance $\td{\ve}_0$ is compact supported, then finite speed of propagation implies that the support of $\td{\ve}$ will keep compact during the lifespan, which provides crucial convenience when performing energy estimate. At last, we apply the weight Sobolev space and the bootstrap framework established in \cite{Collot2018} to prove the main theorem.

Compared with earlier works on formation, persistence and continuation of shocks of conservation laws, this article mainly stresses emphasis on the precise descriptions on the shape and dynamical behavior. Besides, the initial data is assumed to be smooth and compactly supported, which is different from the classic Riemann problem. In particular, the solution we construct is stable under a small perturbation on initial data restricted in a compact set. 

\textbf{Notation} In this article, we use ``$A\les B$'' and ``$A\gtr B$'' to denote $A\leq CB$ and $A\geq CB$, where the constant $C$ may change from line to line and is independent of $\de$.

\section{The Evolution of $\ve$ and $Y_l$}

\subsection{The Evolution of $\ve$}

Now we introduce the self-similar variables corresponding to Burgers equation.

First we let $t\in[T-\de,T)$ and $s=-\log(T-t)$, so $s\in[-\log\de,+\infty).$ Then we define
$$\al_l:=\f{2i_l+1}{2i_l}=1+\f{1}{2i_l},$$
then $1<\al_l<2$. We also define
$$\bar{X}_l:=\f{x}{(T-t)^{\al_l}},\ \ \ X:=\f{x}{(T-t)^{\al_1}},$$
$$\bar{Y_l}:=\f{y_l(t)}{(T-t)^{\al_l}},\ \ \ Y_l:=\f{y_l(t)}{(T-t)^{\al_1}},$$
then 
\begin{equation}\label{trans}
\bar{X}_l=X\f{1}{(T-t)^{\al_l-\al_1}}=e^{(\al_l-\al_1)s}X,\ \ \ \p_X=e^{(\al_l-\al_1)s}\p_{\bar{X}_l}.
\end{equation}

We define $\ve$ by the relation
\begin{equation} \label{ve}
\ve(s,X):=e^{\al_1s}\td{\ve}(s, X),
\end{equation} then direct computation shows that
\begin{align*}
\p_tu_l+u_l\p_xu_l=&(T-t)^{(\f{1}{2i_l}-1)}\Bigg[-\f{1}{2i_l}\Psi_{i_l}(\bar{X}_l-\bar{Y}_l)+\al_l(\bar{X}_l-\bar{Y}_l)\p_{\bar{X}_l}\Psi_{i_l}(\bar{X}_l-\bar{Y}_l)\\
&+\Psi_{i_l}(\bar{X}_l-\bar{Y}_l)\p_{\bar{X}_l}\Psi_{i_l}(\bar{X}_l-\bar{Y}_l)\Bigg]-\f{1}{T-t}\f{dy_l}{dt} \p_{\bar{X}_l}\Psi_{i_l}(\bar{X}_l-\bar{Y}_l).
\end{align*}
Note that since $Y_l=\f{y_l(t)}{(T-t)^{\al_1}}$ and $y_l=e^{-\al_1 s}Y_l$, we have
\begin{equation} \label{yt}
\f{dy_l}{dt}=\f{dy_l}{ds}\f{ds}{dt}=\f{dy_l}{ds}\f{1}{T-t}=e^s(\f{dY_l}{ds}e^{-\al_1s}-\al_1e^{-\al_1s}Y_l)=e^{-\f{s}{2i_1}}(\f{dY_l}{ds}-\al_1Y_l).
\end{equation}
Combinning with (\ref{sebur}), (\ref{deco1}), (\ref{ve}) and (\ref{yt}), we get that
\begin{align} \notag
&\p_s \ve-\f{1}{2i_1}\ve+\al_1X\p_X\ve+\ve\p_X\ve+\p_X\ve\sum_{l=1}^{L}e^{(\al_1-\al_l)s}\Psi_{i_l}(\bar{X}_l-\bar{Y}_l)+\ve\sum_{l=1}^{L}\p_{\bar{X}_l}\Psi_{i_l}(\bar{X}_l-\bar{Y}_l)\\ 
&+\sum_{l_1\neq l_2}e^{(\al_1-\al_{l_1})s}\Psi_{i_{l_1}}(\bar{X}_{l_1}-\bar{Y}_{l_1})\p_{\bar{X}_{l_2}}\Psi_{i_{l_2}}(\bar{X}_{l_2}-\bar{Y}_{l_2})
+\sum_{l=1}^{L}(\f{d}{ds}Y_l-\al_1 Y_l)\p_{\bar{X}_l}\Psi_{i_l}(\bar{X}_l-\bar{Y}_l)=0. \label{evo}
\end{align}

\subsection{The Evolution of $Y_l$}
Let's go back to $y_l(t)$, which measures the distance between the centers of self-similar parts $u_1$ and $u_l$. Consider two characteristic lines $x_1(t)$ and $x_l(t)$ starting from $x_1(T-\de)=0$ and $x_l(T-\de)=y_l(T-\de)=y_{l,0},$ then 
$$y_l(t)=x_l(t)-x_1(t)$$ and 
$$\f{dy_l}{dt}=\f{dx_l}{dt}-\f{dx_1}{dt}.$$
On the other hand, by using the characteristic line of Burgers equation,
$$\f{dx_l}{dt}=u(T-\de,x_l(T-\de))=u(T-\de,y_l(T-\de)=u_0(y_{l,0}),$$
$$\f{dx_1}{dt}=u(T-\de,0)=u_0(0),$$
which implies that $h_l:=\f{dy_l}{dt}=u_0(y_{l,0})-u_0(0)$ is a constant. Together with (\ref{yt}), it follows that 
\begin{equation}\label{hl}
e^{-\f{s}{2i_1}}(\f{dY_l}{ds}-\al_1Y_l)=h_l.
\end{equation}
Integrating it from $s=-\log\de$, we get 
\begin{equation}\label{Yls}
Y_l(s)=e^{\al_1s}\de^{\al_1}Y_{l,0}+e^{\al_1s}(\de-e^{-s})h_l,
\end{equation}
where $Y_{l,0}=Y_l(-\log\de)=y_{l,0}\de^{-\al_1}=3(l-1)\de^{-1}\geq0$.

About $h_l$ we have the following lemma.
\begin{lem}For any $l=2,...,L-1$, we have $h_l<0$ and $0>h_l>h_{l+1}$  
\begin{equation}\label{hl1}
|h_l|\les\de^{\f{1}{2i_m(2i_m+1)}},
\end{equation}
 where $i_m=\max\{i_1,....,i_L\}.$
\end{lem}
\proof
In order to show $0>h_l>h_{l+1}$, it suffices to show that $\sum_{l=1}^{L}u_{l,0}$ is a strictly decreasing function in $x$. Indeed, for any $i$, $\Psi_{i_l}$ is strictly decreasing, so are $u_{l,0}=\de^{\f{1}{2i_l}}\Psi((x-y_{l,0})\de^{-\al_l})$ and $\sum_{l=1}^{L}u_{l,0}.$ 

By the definition of $h_l$, we have
\begin{align*}
h_l=&u_0(y_{l,0})-u_0(0)\\
=&\sum_{l'=1}^{L}u_{l',0}(y_{l,0})+\td{\ve}_0(y_{l,0})-\sum_{l'=1}^{L}u_{l',0}(0)-\td{\ve}_0(0)\\
=&\sum_{l'=1}^{L}u_{l',0}(y_{l,0})-\sum_{l'=1}^{L}u_{l',0}(0)\\
=&\sum_{l'=1}^{L}(u_{l',0}(y_{l,0})-u_{l',0}(0))\\
<&0,
\end{align*}
where the fact that $y_{l,0}$, $0\notin \mathrm{supp}\ \tilde{\ve}_0$ is used in the third equality.

At last, 
\begin{align*}
h_l=&\sum_{l'=1}^{L}(u_{l',0}(y_{l,0})-u_{l',0}(0))\\
=&\sum_{l'=1}^{L}\de^{\f{1}{2i_{l'}}}\Psi_{i_{l'}}(\f{y_{l,0}-y_{l',0}}{\de^{1+\f{1}{2i_{l'}}}})-\sum_{l'=1}^{L}\de^{\f{1}{2i_{l'}}}\Psi_{i_{l'}}(\f{-y_{l',0}}{\de^{1+\f{1}{2i_{l'}}}})\\
=:&I_1+I_2.
\end{align*}

For $I_2$, as a consequence of  \eqref{yl0} we have
$$|I_2|\les\sum_{l'=1}^{L}\de^{\f{1}{2i_{l'}}}(3(l'-1)\de^{-1})^{\f{1}{2i_{l'}+1}}\les\de^{\f{1}{2i_m(2i_m+1)}}.$$

For $I_1$, we rewrite is as
\begin{align*}
I_1
=&\sum_{i_{l'}<i_l}\de^{\f{1}{2i_{l'}}}\Psi_{i_{l'}}(\f{y_{l,0}-y_{l',0}}{\de^{1+\f{1}{2i_{l'}}}})+\sum_{i_{l'}>i_l}\de^{\f{1}{2i_{l'}}}\Psi_{i_{l'}}(\f{y_{l,0}-y_{l',0}}{\de^{1+\f{1}{2i_{l'}}}}),
\end{align*}
noting that $\f{y_{l,0}-y_{l',0}}{\de^{1+\f{1}{2i_{l'}}}}=3(l-1)\de^{-1+\f{1}{2i_l}-\f{1}{2i_{l'}}}-3(l'-1)\de^{-1}=0$ when $i_{l'}=i_l$.

In the sum $\sum_{i_{l'}<i_l}$, we have $\f{1}{2i_l}-\f{1}{2i_{l'}}<0$, so $\de^{-1+\f{1}{2i_l}-\f{1}{2i_{l'}}}\gg\de^{-1}$. Therefore when $\de$ is small enough, there holds
\begin{align*}
|\sum_{i_{l'}<i_l}\de^{\f{1}{2i_{l'}}}\Psi_{i_{l'}}(3(l-1)\de^{-1+\f{1}{2i_l}-\f{1}{2i_{l'}}}-3(l'-1)\de^{-1})|\les&\sum_{i_{l'}<i_l}(\de^{-1+\f{1}{2i_l}-\f{1}{2i_{l'}}})^{\f{1}{2i_{l'}+1}}\de^{\f{1}{2i_{l'}}}\\
\les&\sum_{i_{l'}<i_l}\de^{\f{1}{(2i_{l'}+1)2i_l}}\\
\les&\de^{\f{1}{(2i_{m}+1)2i_m}}.
\end{align*}

In the sum $\sum_{i_{l'}>i_l}$, we have $\f{1}{2i_l}-\f{1}{2i_{l'}}>0$ so $\de^{-1+\f{1}{2i_l}-\f{1}{2i_{l'}}}<\de^{-1}.$ Like above when $\de$ is small enough, we have
\begin{align*}
|\sum_{i_{l'}>i_l}\de^{\f{1}{2i_{l'}}}\Psi_{i_{l'}}(3(l-1)\de^{-1+\f{1}{2i_l}-\f{1}{2i_{l'}}}-3(l'-1)\de^{-1})|\les&\sum_{i_{l'}>i_l}\de^{\f{1}{2i_{l'}}}(\de^{-1})^{\f{1}{2i_{l'}+1}}\\
\les&\de^{\f{1}{(2i_{m}+1)2i_m}}.
\end{align*}

In conclusion, $|h_l|\les\de^{\f{1}{(2i_{m}+1)2i_m}}$.

\qed

By the above lemma, we get the following estimate
\begin{equation}\label{absorb}
e^{\al_1s}(\de-e^{-s})h_l<\f{1}{2}e^{\al_1s}\de^{\f{1}{2i_1}},
\end{equation}
since $(\de-e^{-s})h_l<\de<\f{1}{4}\de^{\f{1}{2i_1}}$ when $\de$ is small enough and $s>-\log\de$.

\section{Some Useful Facts}
In this section we will record some properties of $\Psi_i$ and $\phi_i$, summarized or proved by Collot, Ghoul and Masmoudi in \cite{Collot2018}.
\begin{prop}
We have the asymptotic property of $\Psi_i(X)$,
$$\Psi_i(X)=-X+X^{2i+1}+\sum_{k=2}^{+\infty}c_{i,l}X^{2ki+1}\ \ as\ \ X\to0$$and
$$\Psi_i(X)=-sgn(X)|X|^{\f{1}{2i+1}}+sgn(X)\f{|X|^{-1+\f{2}{2i+1}}}{2i+1}+O(|X|^{-2+\f{3}{2i+1}})\ \ as\ \ |X|\to+\infty.$$
Moreover, we have 
\begin{equation}\label{psibd}
|\Psi_i(X)|\approx|X|(1+|X|)^{\f{1}{2i+1}-1},\ \ |\p_X\Psi_i(X)|\les|X|^{\f{1}{2i+1}-1},\ \ |\p^2_X\Psi_i(X)|\les|X|^{\f{1}{2i+1}-2}.
\end{equation}
\end{prop}

Now we define the linear operators $H_X$ and $\td{H}_X$, such that
$$H_Xf=-\f{1}{2i_1}f+\al_1X\p_X f+f\p_X\Psi_{i_1}(X)+\Psi_{i_1}(X)\p_Xf$$ and
\begin{align}
\td{H}_Xf=&-\f{1}{2i_1}f+\al_1X\p_X f+\p_Xf\sum_{l=1}^{L}e^{(\al_1-\al_l)s}\Psi_{i_l}(\bar{X}_l-\bar{Y}_l)+f\sum_{l=1}^{L}\p_{\bar{X}_l}\Psi_{i_l}(\bar{X}_l-\bar{Y}_l)\notag\\ \label{tdHX}
=&H_Xf+\p_Xf\sum_{l=2}^{L}e^{(\al_1-\al_l)s}\Psi_{i_l}(\bar{X}_l-\bar{Y}_l)+f\sum_{l=2}^{L}\p_{\bar{X}_l}\Psi_{i_l}(\bar{X}_l-\bar{Y}_l).
\end{align}

The spectrum properties of $H_X$ is proved in \cite{Collot2018}.
\begin{prop}
The point spectrum of $H_X$ is
$$\varUpsilon=\left\{\f{j-2i_1-1}{2},\ \ j\in\mathbb{N} \right\}.$$
In particular,
\begin{equation}\label{HXspec}
H_X(\phi_j)=\f{j-2i_1-1}{2}\phi_j,
\end{equation}
where
\begin{equation}\label{phi}
\phi_j(X)\approx|X|^j(1+|X|)^{\f{j-2}{3}-j},\ \ \ |\f{\p_X\phi_j}{\phi_j}|\les\f{1}{|X|}.
\end{equation}
\end{prop}

We will use the function space
$$\mathcal{B}_{j,q}=\left\{u: \int_{\mathbb{R}}\f{u^{2q}dX}{\phi^{2q}_j(X)|X|}+\int_{\mathbb{R}}\f{(X\p_Xu)^{2q}dX}{\phi^{2q}_j(X)|X|}<+\infty, \ \ q\in\mathbb{N}^*\right\},$$ 
as well as an embedding lemma, whose proof can be found in \cite{Collot2018}.
\begin{lem}
Let $q\in\mathbb{N}^*$. Then for any $u\in \mathcal{B}_{j,q}$ one has
\begin{equation}\label{ebd}
\left\|\f{u}{\phi_j}\right\|_{L^\infty(\mathbb{R})}^{2q}\leq C(q)\left( \int_{\mathbb{R}}\f{u^{2q}dX}{\phi^{2q}_j(X)|X|}+\int_{\mathbb{R}}\f{(X\p_Xu)^{2q}dX}{\phi^{2q}_j(X)|X|}\right).
\end{equation}
\end{lem}

\section{Bootstrap Argument}
Since all the preparation of bootstrap have been made, we now state main theorem precisely.
\begin{thm}
Let $j=2i_1+2$, 
and there exist $0<\de, \ka\ll1$, $q\gg1$ such that for all
 $\ve_0(X)$ supported in $X\in[\f{4}{3\de},\f{5}{3\de}]$ with 
 \begin{align} \label{ininorm}
 \|\ve_0(X)\|_{\mathcal{B}_{2i_1+2,q}}\leq\de^{(\f{1}{2}-\ka)},
 \end{align}
  then the solution to equation \emph{(\ref{evo})} with data $\ve_0$ is global, where $\ve$ is defined in \emph{(\ref{ve})}.
 
In particular, if we take 
$$h_l=\sum_{l'=1}^{L}\de^{\f{1}{2i_{l'}}}\left(\Psi_{i_{l'}}\left(\f{y_{l,0}-y_{l',0}}{\de^{1+\f{1}{2i_{l'}}}}\right)-\Psi_{i_{l'}}\left(\f{-y_{l',0}}{\de^{1+\f{1}{2i_{l'}}}}\right)\right)$$ with 
$y_{l,0}=3(l-1)\de^{\f{1}{2i_1}},$
 then the equation \emph{(\ref{bur})} admits solutions of the form 
$$u(t,x)=\sum_{l=1}^{L}(T-t)^{\f{1}{2i_l}}\Psi_{i_l}(\f{x-y_l(t)}{(T-t)^{\al_l}})+\td{\ve}(t,x),$$
where $y_l(t)$ satisfies
$$y_l(t)=3(l-1)\de^{\f{1}{2i_1}}+(\de-(T-t))\de^{\f{1}{2i_1}}h_l,$$
and $\Psi_{i_l}(\bar{X}_l)$ is the solution of $-\f{1}{2i_l}\Psi_{i_l}+\al_l\bar{X}_l\p_{\bar{X}_l}\Psi_{i_l}+\Psi_{i_l}\p_{\bar{X}_l}\Psi_{i_l}=0$ with $\bar{X}_l=e^{(\al_l-\al_1)s}X$.
\end{thm}

Next we need to estimate the support and weight Sobolev norm of $\ve$ by bootstrap argument. 

\begin{prop}(The support of $\ve(s,X)$)For any $s>-\log\de$, we have
\begin{equation}\label{sptX}
\emph{supp}\{\ve(s,X)\}\subseteq\mathcal{X}(s)=\left\{X: e^{s}\leq X \leq2\de^{\f{1}{2i_1}}e^{\al_1s}  \right\},
\end{equation}
so that by \emph{(\ref{Yls})} and (\ref{absorb})
\begin{equation}\label{X-Y}
|X-Y_l|\approx
\begin{cases}
 \de^{\f{1}{2i_1}}e^{\al_1s},\ \ l\geq2,\\
 e^s,\ \ l=1.
\end{cases}
\end{equation}
\end{prop}

Next we define the vector field $A$,
\begin{align*}
A=&\al_1X\p_X+\sum_{l=1}^{L}e^{(\al_1-\al_l)s}\Psi_{i_l}(\bar{X}_l-\bar{Y}_l)\p_X\\
=&\al_1X\p_X+\sum_{l=1}^{L}e^{(\al_1-\al_l)s}\Psi_{i_l}(e^{(\al_l-\al_1)s}({X}-{Y}_l))\p_X.
\end{align*}

\begin{prop}(The $\mathcal{B}_{2i_1+2,q}$ norm of $\ve(s,X)$)
There exist constants $q\in\mathbb{N}^*$, $K_0,$ $K_1\gg1$ and $\de,\ \ka\ll1$ such that if the solution $\ve$ to (\ref{evo}) with initial data $\ve_0=\ve(-\log\de)$  satisfying
\begin{equation}
\int_{\mathbb{R}}\f{\ve_0^{2q}+(A\ve_0)^{2q} }{\phi_{2i_1+2}^{2q}(X)|X| }dX\leq e^{-{2q(\f{1}{2}-\ka)s_0}},
\end{equation}
then the solution $\ve$ is global and satisfies
\begin{equation}\label{btbd1}
\left(\int_{\mathbb{R}}\f{\ve^{2q}}{\phi^{2q}_{2i_1+2}(X)|X|}dX  \right)^{\f{1}{2q}}\leq K_0e^{-(\f{1}{2}-\ka)s}
\end{equation}
and
\begin{equation}\label{btbd2}
\left(\int_{\mathbb{R}}\f{(A\ve)^{2q}}{\phi^{2q}_{2i_1+2}(X)|X|}dX  \right)^{\f{1}{2q}}\leq K_1e^{-(\f{1}{2}-\ka)s}.
\end{equation}
\end{prop}

\begin{cor}
For $s\geq-\log\de$ and $X\in\mathcal{X}(s),$ with the ``$\les$'' depending on bootstrap constants $K_0$ and $K_1$, there hold
\begin{equation}\label{vebd}
|\ve|\les|X|e^{-(\f{1}{2}-\ka)s}
\end{equation}
and
\begin{equation}\label{pvebd}
|\p_X\ve|\les e^{-(\f{1}{2}-\ka)s}.
\end{equation}
\end{cor}

\subsection{The Evolution of the Support}
Before we prove \textbf{Proposition 4.2}, an estimate measuring the lower bound of transport speed of (\ref{evo}) is needed.

\begin{lem}For any $s\geq-\log\de$ with $\de$ small enough and $$X\in\mathcal{X}(s)=
\{e^s\leq X \leq 2 \de^{\f{1}{2i_1}}e^{\al_1s} \},$$
we have
\begin{equation}\label{lowbd}
X\leq\al_1X+\sum_{l=1}^{L}e^{(\al_1-\al_l)s}\Psi_{i_l}(e^{(\al_l-\al_1)s}({X}-{Y}_l(s))).
\end{equation}
\end{lem}
\proof
First of all, note that when $X\in\mathcal{X}(s)$, the term when indexed by $l=1$ is the unique negative term in the sum, since $\Psi_{i_l}$ is decreasing odd function in $X$. So it suffices to prove
$$\Psi_{i_1}(X)>-\f{1}{2i_1}X,$$ where we note that $Y_1(s)=0.$

Then for fixed $s$, it suffices to show that
\begin{equation*}
\inf_{X\in\mathcal{X}(s)}\left\{\Psi_{i_1}(X)\right\}>\sup_{X\in\mathcal{X}(s)}\left\{-\f{1}{2i_1}X\right\}=-\f{1}{2i_1}e^s.
\end{equation*}

Indeed by (\ref{sptX}), (\ref{X-Y}) and (\ref{psibd}), we have 
\begin{align*}
\inf_{X\in\mathcal{X}(s)}\left\{\Psi_{i_1}(X)\right\}=&
\Psi_{i_1}(2\de^{\f{1}{2i_1}}e^{\al_1s}))\\
\geq&-c(\de^{\f{1}{2i_1}}e^{\al_1s})^{\f{1}{2i_1+1}}\\
\geq&-c\de^{\f{1}{2i_1(2i_1+1)}}e^{(\al_1-1)s}\\
\geq&-\f{1}{2i_1}e^s,
\end{align*}
where $\de$ is small enough.

\qed

\proof[Proof of Propsition 4.2]

Under the self-similar variables, by (\ref{sptx}), the support of $\ve(-\log\de,X_0)$ is 
\begin{equation}\label{sptX0}
\mathcal{X}_0=\left\{X_0:\f{4}{3\de}\leq X_0\leq\f{5}{3\de}\right\}.
\end{equation}
Next we use Lagrangian trajectories $\Phi^{X_0}(s)$ defined by
\begin{equation}\label{phiX}
\begin{cases}
\f{d}{ds}\Phi^{X_0}(s)=\al_1\Phi^{X_0}(s)+\sum_{l=1}^{L}e^{(\al_1-\al_l)s}\Psi_{i_l}\left(e^{(\al_l-\al_1)s}(\Phi^{X_0}(s)-Y_l(s))\right)+\ve,\\
\Phi^{X_0}(-\log\de)=X_0.
\end{cases}
\end{equation}

For the lower bound,  by (\ref{lowbd}) (\ref{sptX}) and (\ref{vebd}), we have
\begin{align*}
\f{d}{ds}\Phi^{X_0}(s)\geq&\Phi^{X_0}(s)-|\Phi^{X_0}(s)||\ve|\\
\geq&\Phi^{X_0}(s)-2C_K\de^{\f{1}{2i_1}}e^{(\al_1-\f{1}{2}+\ka)s}.
\end{align*}
Integrating from $-\log\de,$ we have 
\begin{align*}
\Phi^{X_0}(s)\geq&|X_0|\de e^{s}-2C_K\de^{\f{1}{2i_1}}e^{s}\int_{-\log\de}^{s}e^{(\al_1-\f{3}{2}+\ka)\td{s}}d\td{s}\\
\geq&\f{4}{3}e^{s}-2C_K\de^{\f{1}{2i_1}}e^{s}\de^{-(\al_1-\f{3}{2}+\ka)}\\
\geq&(\f{4}{3}-2C_K\de^{\f{1}{2}-\ka})e^s\\
\geq&e^s,
\end{align*}
as $\de$ is small enough.

For the upper bound, by (\ref{sptX}) (\ref{Yls}) and (\ref{psibd}) , we have
\begin{align}
|e^{(\al_1-\al_l)s}\Psi_{i_l}\left(e^{(\al_l-\al_1)s}(\Phi^{X_0}(s)-Y_l(s))\right)|\leq& |e^{(\al_1-\al_l)s}\Psi_{i_l}(e^{(\al_l-\al_1)s}\de^{\f{1}{2i_1}}e^{\al_1s})|\notag\\
\leq&|Ce^{(\al_1-\al_l)s}\Psi_{i_l}(\de^{\f{1}{2i_l}}e^{\al_ls})|\notag\\
\leq&Ce^{(\al_1-\al_l)s}\de^{\f{1}{2i_l(2i_l+1)}}e^{(\al_l-1)s}\notag\\
\leq&C\de^{\f{1}{2i_m(2i_m+1)}}e^{(\al_1-1)s}\label{psil},
\end{align}
where $i_m=\max\{i_1,...,i_L \}.$
Then (\ref{phiX}) and (\ref{vebd}) implies that 
\begin{align*}
\f{d}{ds}\Phi^{X_0}(s)\leq \al_1\Phi^{X_0}(s) +LC\de^{\f{1}{2i_m(2i_m+1)}}e^{(\al_1-1)s}+2C_K\de^{\f{1}{2i_1}}e^{(\al_1-\f{1}{2}+\ka)s}.
\end{align*}
Integrate it then we get
\begin{align*}
\Phi^{X_0}(s)\leq&\f{5}{3}\de^{\f{1}{2i_1}}e^{\al_1s}+LC\de^{1+\f{1}{2i_m(2i_m+1)}}e^{\al_1s}+2C_K\de^{\f{1}{2i_1}+\f{1}{2}-\ka}e^{\al_1s}\\
\leq&2\de^{\f{1}{2i_1}}e^{\al_1s},
\end{align*}
as $\de$ is small enough.

\qed

\proof[Proof of Corollary 4.4]
By (\ref{lowbd}) we deduce
$$|\ve|+|X\p_X\ve|\leq|\ve|+|A\ve|,$$
therefore the bootstrap assumption (\ref{btbd1}) and (\ref{btbd2}) imply that 
$$\int_{\mathbb{R}}\f{ \ve^{2q}  }{\phi^{2q}_{2i_1+2}(X)|X|}+\int_{\mathbb{R}}\f{ (X\p_X\ve)^{2q}  }{\phi^{2q}_{2i_1+2}(X)|X|}\les e^{-(\f{1}{2}-\ka)s}, $$
According to the embedding lemma (\ref{ebd}) we have $|\ve|\les e^{-(\f{1}{2}-\ka)s}|\phi_{2i_1+2}|\les e^{-(\f{1}{2}-\ka)s}|X|.$ Similarly, $|\p_X\ve|\les e^{-(\f{1}{2}-\ka)s}|\phi_{2i_1+2}||X|^{-1}\les e^{-(\f{1}{2}-\ka)s}.$

\qed

\subsection{Energy Estimate}
The proof of \textbf{Proposition 4.3} is split into several lemmas.
\begin{lem}
Assume that $u$ and $\Theta$ are smooth and satisfy
$$u_s+\td{H}_Xu=\Theta$$for $s\in[-\log\de,+\infty),$ and $q\in\mathbb{N}^*$ large enough, where $\td{H}_X$ is defined in \emph{(\ref{tdHX})}. Then there exist $C>0$ independent of q such that for $\de$ small enough the following energy estimate holds,
\begin{equation}\label{engy}
\f{d}{ds}\left( \int_{\mathbb{R}}\f{ u^{2q} }{\phi^{2q}_{2i_1+2}(X)|X|}dX \right)\leq -\left(\f{1}{2}-\f{C}{q}   \right)\int_{\mathbb{R}}\f{ u^{2q} }{\phi^{2q}_{2i_1+2}(X)|X|}dX+\int_{\mathbb{R}}\f{ u^{2q-1}\Theta }{\phi^{2q}_{2i_1+2}(X)|X|}dX
\end{equation}
\end{lem}
\proof
In this proof, $0<\mu\ll1$, which may change line by line. By (\ref{evo})
\begin{align*}
\f{d}{ds}\left( \f{1}{2q}\int_{\mathbb{R}}\f{ u^{2q} }{\phi^{2q}_{2i_1+2}(X)|X|}dX  \right)=&\int_{\mathbb{R}}\f{ u^{2q-1}u_s }{\phi^{2q}_{2i_1+2}(X)|X|}dX\\
=&\int_{\mathbb{R}}\f{ u^{2q-1} }{\phi^{2q}_{2i_1+2}(X)|X|}\Big[
-H_Xu-u\sum_{l=2}^{L}\p_{\bar{X}_l}\Psi_{i_l}(\bar{X}_l-\bar{Y}_l)\\
&-\p_Xu\sum_{l=2}^{L}e^{(\al_1-\al_l)s}\Psi_{i_l}(\bar{X}_l-\bar{Y}_l)+\Theta \Big] dX\\
 =&I_1+I_2+I_3+I_4.
\end{align*}

For $I_1$, by (\ref{HXspec}) and the similar lemma in \cite{Collot2018}, we get $I_1\leq-\left( \f{1}{2}-\f{C}{q} \right)\int_{\mathbb{R}}\f{ u^{2q} }{\phi^{2q}_{2i_1+2}(X)|X|}dX.$

For $I_2$, by (\ref{trans}), (\ref{X-Y}) and (\ref{psibd}), we have
\begin{equation}\label{pxpsi}
\p_{\bar{X}_l}\Psi_{i_l}(\bar{X}_l-\bar{Y}_l)\les(e^{(\al_l-\al_1) s}e^s)^{\f{1}{2i_l+1}-1}\les e^{(\al_l-\al_1+1)(\f{1}{2i_l+1}-1)s}\les e^{(\f{\al_1}{\al_l}-1-\f{1}{\al_l})s}\les \de^{\mu},
\end{equation}
where $\f{\al_1}{\al_l}-1-\f{1}{\al_l}<0$ and $\mu>0$ is a proper chosen constant. So 
\begin{align*}
|I_2|\les \de^{\mu}\int_{\mathbb{R}}\f{ u^{2q} }{\phi^{2q}_{2i_1+2}(X)|X|}dX.
\end{align*}

For $I_3$, by integration by parts we note that 
\begin{align*}
&\int_{\mathbb{R}}\f{ u^{2q-1} \p_Xu}{\phi^{2q}_{2i_1+2}(X)|X|} e^{(\al_1-\al_l)s}\Psi_{i_l}(\bar{X}_l-\bar{Y}_l) dX\\=&
\f{1}{2q}\int_{\mathbb{R}}\f{ \p_X(u^{2q}) }{\phi^{2q}_{2i_1+2}(X)|X|}e^{(\al_1-\al_l)s}\Psi_{i_l}(\bar{X}_l-\bar{Y}_l)dX\\
=&\f{1}{2q}\int_{\mathbb{R}}\f{ u^{2q}(\p_X\phi_{2i_1+2}) }{\phi^{2q+1}_{2i_1+2}(X)|X|}e^{(\al_1-\al_l)s}\Psi_{i_l}(\bar{X}_l-\bar{Y}_l)dX-\f{1}{2q}\int_{\mathbb{R}}\f{ u^{2q} }{\phi^{2q}_{2i_1+2}(X)}\p_X \left( \f{e^{(\al_1-\al_l)s}\Psi_{i_l}(\bar{X}_l-\bar{Y}_l)}{|X|}   \right) dX\\
=&I_{3,1}+I_{3,2}.
\end{align*}

For each $l$ in the sum of the $I_{3,1}$ we use (\ref{trans}), (\ref{phi}) and get
$$e^{(\al_1-\al_l)s}\f{ \Psi_{i_l}(\bar{X}_l-\bar{Y}_l)}{|X|}\les e^{(\f{\al_1}{\al_l}-1)s}\f{|X-Y_l|^{\f{1}{2i_l+1}}}{|X|}.$$
 We estimate for two cases:

\emph{Case1.} $|X|\gg|Y|$ or $|X|\approx|Y|$. At this time by (\ref{sptX}) and (\ref{X-Y}) we have
$$e^{(\f{\al_1}{\al_l}-1)s}\f{|X-Y_l|^{\f{1}{2i_l+1}}}{|X|}\les e^{(\f{\al_1}{\al_l}-1)s}|X|^{-\f{1}{\al_l}}\les e^{(\f{\al_1-1}{\al_l}-1)s}\les \de^{\mu}.
$$

\emph{Case2} .$|X|\ll|Y|$. At this time by (\ref{sptX}) and (\ref{X-Y}) we have
$$e^{(\f{\al_1}{\al_l}-1)s}\f{|X-Y_l|^{\f{1}{2i_l+1}}}{|X|}\les e^{(\f{\al_1}{\al_l}-1)s}\f{ (\de^{\f{1}{2i_1}}e^{\al_1s})^{\f{1}{2i_l+1}} }{e^s}\les \de^{\f{1}{2i_1(2i_l+1)}}e^{(\al_1-2)s}\les \de^{\mu}.$$

So by (\ref{phi}), $$I_{3,1}\les\de^{\mu}\int_{\mathbb{R}}\f{ u^{2q} }{\phi^{2q}_{2i_1+2}(X)|X|}dX.$$

The estimate of $I_{3,2}$ is similar to those of $I_2$ and $I_{3,1}$:
$$\p_X \left( \f{e^{(\al_1-\al_l)s}\Psi_{i_l}(\bar{X}_l-\bar{Y}_l)}{|X|}   \right)=\f{\p_{\bar{X}_l}\Psi_{i_l}(\bar{X}_l-\bar{Y}_l)}{|X|}
-\f{e^{(\al_1-\al_l)s}\Psi_{i_l}(\bar{X}_l-\bar{Y}_l)}{|X|^2}\les\de^{\mu},
$$
Thus 
\begin{align*}
|I_3|\les \de^{\mu}\int_{\mathbb{R}}\f{ u^{2q} }{\phi^{2q}_{2i_1+2}(X)|X|}dX.
\end{align*}

$I_4$ is trivial.Now we choose proper $q$, $C$ and obtain the (\ref{engy}).

\qed

\begin{lem}
There exists $K_0\gg1$ large enough independent of $K_1$ such that for $s>s_0$ we have
$$
\left(\int_{\mathbb{R}}\f{ \ve^{2q} }{\phi^{2q}_{2i_1+2}(X)|X|}dX  \right)^{\f{1}{2q}}\leq\f{K_0}{2}e^{-(\f{1}{2}-\ka)s}.
$$
\end{lem}
\proof In this proof, $$\td{\al}:=\f{1}{2i_1}+(\f{\al_1}{\al_l}-1-\f{1}{\al_l})<0$$ and $C$ is a constant which may change from line to line. 
We apply (\ref{engy}) to (\ref{evo}) and get
\begin{align*}
&\f{d}{ds}\left(\f{1}{2q} \int_{\mathbb{R}}\f{ \ve^{2q} }{\phi^{2q}_{2i_1+2}(X)|X|}dX  \right)\\
\leq &\left( \f{1}{2}-\f{C}{q}  \right)\int_{\mathbb{R}}\f{ \ve^{2q} }{\phi^{2q}_{2i_1+2}(X)|X|}dX -\int_{\mathbb{R}}\f{\ve^{2q-1}}{\phi^{2q}_4(X)|X|}dX 
\Big[\sum_{l=1}^{L}(\f{d}{ds}Y_l-\al_1 Y_l)\p_{\bar{X}_l}\Psi_{i_l}(\bar{X}_l-\bar{Y}_l)\\
&+\sum_{l_1\neq l_2}e^{(\al_1-\al_{l_1})s}\Psi_{i_{l_1}}(\bar{X}_{l_1}-\bar{Y}_{l_1})\p_{\bar{X}_{l_2}}\Psi_{i_{l_2}}(\bar{X}_{l_2}-\bar{Y}_{l_2})+\ve\p_X\ve\Big]\\
=&\left( \f{1}{2}-\f{C}{q}  \right)\int_{\mathbb{R}}\f{ \ve^{2q} }{\phi^{2q}_{2i_1+2}(X)|X|}dX-(I_1+I_2+I_3).
\end{align*}

For $I_1$, due to (\ref{hl}) and (\ref{hl1}), we have $|\f{d}{ds}Y_l-\al_1 Y_l|=|e^{\f{1}{2i_1}s}h_l|\les e^{\f{1}{2i_1}s}\de^{\f{1}{2i_m(2i_m+1)}}$. Then by (\ref{pxpsi}) and Young's inequality,
\begin{align*}
|I_1|\les& 
\left| \int_{\mathbb{R}}\f{\ve^{2q-1}}{\phi^{2q}_{2i_1+2}(X)|X|}
e^{\f{1}{2i_1}s}\de^{\f{1}{2i_m(2i_m+1)}} \sum_{l=1}^{L} \p_{\bar{X}_l}\Psi_{i_l}(\bar{X}_l-\bar{Y}_l)dX\right|\\
 \les&\f{c_0}{q} \int_{\mathbb{R}}\f{\ve^{2q}}{\phi^{2q}_{2i_1+2}(X)|X|}dX+\sum_{l=1}^{L}C(q)\int_{\mathcal{X}(s)}\f{ (e^{\f{1}{2i_1}s}\de^{\f{1}{2i_m(2i_m+1)}}   e^{(\f{\al_1}{\al_l}-1-\f{1}{\al_l})s}   )^{2q} }{\phi^{2q}_{2i_1+2}(X)|X| }dX, 
\end{align*}
where $c_0$ is a small number. Next we integrate the second term
\begin{align*}
&\sum_{l=1}^{L}C(q)\int_{\mathcal{X}(s)}\f{ (e^{\f{1}{2i_1}s}\de^{\f{1}{2i_m(2i_m+1)}}e^{(\f{\al_1}{\al_l}-1-\f{1}{\al_l})s})^{2q}  }{\phi^{2q}_{2i_1+2}(X)|X| }dX\\
\les&\sum_{l=1}^{L}(e^{\f{1}{2i_1}s}\de^{\f{1}{2i_m(2i_m+1)}}e^{(\f{\al_1}{\al_l}-1-\f{1}{\al_l})s})^{2q} \de^{\f{1}{2i_1}}e^{\al_1s}\\
\les&\de^{\f{q}{i_m(2i_m+1)}+\f{1}{2i_1}}e^{2q\td{\al}s+\al_1s},
\end{align*}
where $q\gg1.$ 

For $I_2$, by (\ref{psil}) and (\ref{pxpsi}),
\begin{align*}
e^{(\al_1-\al_{l_1})s}\Psi_{i_{l_1}}(\bar{X}_{l_1}-\bar{Y}_{l_1})\p_{\bar{X}_{l_2}}\Psi_{i_{l_2}}(\bar{X}_{l_2}-\bar{Y}_{l_2})\les
\de^{\f{1}{2i_m(2i_m+1)}}e^{(\al_1-1)s}e^{(\f{\al_1}{\al_{l_2}}-1-\f{1}{\al_{l_2}})s},
\end{align*}
so by Young's inequality and (\ref{sptX}), we have
\begin{align*}
|I_2|\les&\f{c_0}{q}\int_{\mathbb{R}}\f{\ve^{2q}}{\phi^{2q}_{2i_1+2}(X)|X|}dX+\sum_{l_1\neq l_2}C(q)  \int_{\mathcal{X}(s)}\f{( \de^{\f{1}{2i_m(2i_m+1)}}e^{(\al_1-1)s}e^{(\f{\al_1}{\al_{l_2}}-1-\f{1}{\al_{l_2}})s})^{2q}}{\phi^{2q}_{2i_1+2}(X)|X|}dX\\
\les&\f{c_0}{q}\int_{\mathbb{R}}\f{\ve^{2q}}{\phi^{2q}_{2i_1+2}(X)|X|}dX+\de^{\f{1}{i_m(2i_m+1)}+\f{1}{2i_1}}e^{2q\td{\al}s+\al_1s}.
\end{align*}

For $I_3$, by (\ref{pvebd}),
$$|I_3|\les\left| \int_{\mathbb{R}}\f{\ve^{2q}\p_X\ve}{\phi^{2q}_{2i_1+2}(X)|X|}dX\right|\les e^{-(\f{1}{2}-\ka)s}\int_{\mathbb{R}}\f{\ve^{2q}}{\phi^{2q}_{2i_1+2}(X)|X|}dX.  $$

At last, 
\begin{align*}
&\f{d}{ds}\left(\f{1}{2q} \int_{\mathbb{R}}\f{ \ve^{2q} }{\phi^{2q}_{2i_1+2}(X)|X|}dX  \right)\\
\leq&\left( \f{1}{2}-\f{C}{q}  \right)\int_{\mathbb{R}}\f{ \ve^{2q} }{\phi^{2q}_{2i_1+2}(X)|X|}dX+(\f{2c_0}{q}+e^{-(\f{1}{2}-\ka)s})\int_{\mathbb{R}}\f{ \ve^{2q} }{\phi^{2q}_{2i_1+2}(X)|X|}dX\\
&+C\de^{\f{1}{i_m(2i_m+1)}+\f{1}{2i_1}}e^{2q\td{\al}s+\al_1s}\\
\leq&\left( \f{1}{2}-\f{C}{q}  \right)\int_{\mathbb{R}}\f{ \ve^{2q} }{\phi^{2q}_{2i_1+2}(X)|X|}dX+C\de^{\f{1}{i_m(2i_m+1)}+\f{1}{2i_1}}e^{2q\td{\al}s+\al_1s}.
\end{align*}

Then take $q$ large enough so that $|C/q|\leq\ka$ and integrate it, we get
\begin{align*}
\int_{\mathbb{R}}\f{ \ve^{2q} }{\phi^{2q}_{2i_1+2}(X)|X|}dX\leq&
e^{-2q(\f{1}{2}-\ka)(s-s_0)}\int_{\mathbb{R}}\f{ \ve_0^{2q} }{\phi^{2q}_{2i_1+2}(X)|X|}dX+Ce^{-2q(\f{1}{2}-\ka)s}\de^{p_1-p_2-(\f{1}{2}-\ka)}\\
\leq&e^{-2q(\f{1}{2}-\ka)s}(1+\de^{p_1-p_2-(\f{1}{2}-\ka)})\\
\leq&e^{-2q(\f{1}{2}-\ka)s}\f{K_0^{2q}}{2^{2q}},
\end{align*}
where $p_1=\f{1}{i_m(2i_m+1)}+\f{1}{2i_1}$, $p_2=2q\td{\al}+\al_1$ and the last inequality holds by taking $K_0$ large enough in terms of $\de$ and independently of $K_1$.

\qed

\begin{lem}
There exists $K_1\gg1$ large enough depending on $K_0$ such that for $s>s_0$ we have
$$\left(\int_{\mathbb{R}}\f{ (A\ve)^{2q} }{\phi^{2q}_{2i_1+2}(X)|X|}dX\right)^{\f{1}{2q}}\leq\f{K_1}{2}e^{-(\f{1}{2}-\ka)s}.$$
\end{lem}
\proof
In this proof, $0<\nu\ll1$ is a small number and
$$\td{\al}:=\max_{l,l_1,l_2}\left\{ (\al_l-\al_1+1)(\f{1}{2i_l+1}-2)+\al_l+\f{1}{2i_1},\ \ (\al_1+\f{\al_1}{\al_{l_1}}+\f{\al_{1}}{\al_{l_2}}-2-\f{1}{\al_{l_1}}-\f{1}{\al_{l_2}})  \right\},$$
whose sign is not important here.

Let $w=:A\ve$, and note that $A$ is the linear transport field of (\ref{evo}), so we have
$$\left[A, \td{H}_X \right]\ve=\left[A, \sum_{l=1}^{L}\p_{\bar{X}_l} \Psi_{i_l}(\bar{X}_l-\bar{Y}_l) \right]\ve=\left(A\sum_{l=1}^{L}\p_{\bar{X}_l} \Psi_{i_l}(\bar{X}_l-\bar{Y}_l)\right)\ve,$$
and get 
\begin{align*}
w_s+\td{H}_Xw+A\sum_{l=1}^{L}\p_{\bar{X}_l} \Psi_{i_l}(\bar{X}_l-\bar{Y}_l)\ve+\ve\p_Xw+w\p_X\ve+\ve[A,\p_X]\ve\\
+A\sum_{l_1\neq l_2}e^{(\al_1-\al_{l_1})s}\Psi_{i_{l_1}}(\bar{X}_{l_1}-\bar{Y}_{l_1})\p_{\bar{X}_{l_2}}\Psi_{i_{l_2}}(\bar{X}_{l_2}-\bar{Y}_{l_2})\\ 
-A\sum_{l=1}^{L}(\f{d}{ds}Y_l-\al_1 Y_l)\p_{\bar{X}_l}\Psi_{i_l}(\bar{X}_l-\bar{Y}_l)=0.
\end{align*}

By energy estimate (\ref{engy}),
\begin{align*}
&\f{d}{ds}\left( \int_{\mathbb{R}}\f{ w^{2q} }{\phi^{2q}_{2i_1+2}(X)|X|}dX \right)\leq -\left(\f{1}{2}-\f{C}{q}   \right)\int_{\mathbb{R}}\f{ w^{2q} }{\phi^{2q}_{2i_1+2}(X)|X|}dX\\+&\int_{\mathbb{R}}\f{ w^{2q-1}}{\phi^{2q}_{2i_1+2}(X)|X|}\Big[   A\sum_{l=1}^{L}\p_{\bar{X}_l} \Psi_{i_l}(\bar{X}_l-\bar{Y}_l)\ve+\ve\p_Xw+w\p_X\ve+\ve[A,\p_X]\ve\\
+&A\sum_{l_1\neq l_2}e^{(\al_1-\al_{l_1})s}\Psi_{i_{l_1}}(\bar{X}_{l_1}-\bar{Y}_{l_1})\p_{\bar{X}_{l_2}}\Psi_{i_{l_2}}(\bar{X}_{l_2}-\bar{Y}_{l_2})\\ 
-&A\sum_{l=1}^{L}(\f{d}{ds}Y_l-\al_1 Y_l)\p_{\bar{X}_l}\Psi_{i_l}(\bar{X}_l-\bar{Y}_l) \Big]dX\\
=& -\left(\f{1}{2}-\f{C}{q}   \right)\int_{\mathbb{R}}\f{ w^{2q} }{\phi^{2q}_{2i_1+2}(X)|X|}dX+(I_1+\cdots+I_6).
\end{align*}

For $I_1$, by (\ref{psil}), (\ref{pxpsi}), (\ref{psibd}) and (\ref{X-Y}) we have
\begin{align}
|A\sum_{l=1}^{L}\p_{\bar{X}_l} \Psi_{i_l}(\bar{X}_l-\bar{Y}_l)|
\les&|\al_1X+\sum_{l=1}^{L}e^{(\al_1-\al_l)s}\Psi_{i_l}(\bar{X}_l-\bar{Y}_l)|\p_X \p_{\bar{X}_l}\Psi_{i_l}(\bar{X}_l-\bar{Y}_l)|\notag\\
\les&|\al_1X+\sum_{l=1}^{L}e^{(\al_1-\al_l)s}\Psi_{i_l}(\bar{X}_l-\bar{Y}_l)|e^{(\al_l-\al_1)s}\p^2_{\bar{X}_l}\Psi_{i_l}(\bar{X}_l-\bar{Y}_l)|\notag\\
\les&(\de^{\f{1}{2i_1}}e^{\al_1s}+\de^{\f{1}{2i_m(2i_m+1)}}e^{(\al_1-1)s})e^{(\al_l-\al_1)s}(e^{s+(\al_l-\al_1)s})^{\f{1}{2i_l+1}-2} \notag\\
\les&\de^{\f{1}{2i_1}}e^{(\al_l-\al_1+1)(\f{1}{2i_l+1}-2)s+\al_ls}\label{Apsi}\\
\les&\de^{\nu}\notag,
\end{align}
where we note that $(\al_l-\al_1+1)(\f{1}{2i_l+1}-2)+\al_l=-2+\f{\al_1}{\al_l}-\f{1}{\al_l}+\al_1=-2+\f{1}{\al_l}(\al_1-1)+\al_1\leq-2+\f{1}{\al_m}\f{1}{2}+\f{3}{2}<0,$ since $1<\al_m\leq\al_l, \al_1\leq\f{3}{2}.$
Therefore, 
\begin{align*}
|I_1|\les &\de^{\nu}\left| \int_{\mathbb{R}}\f{ w^{2q} }{\phi^{2q}_{2i_1+2}(X)|X|}dX\right|^{\f{2q-1}{2q}}\left| \int_{\mathbb{R}}\f{ \ve^{2q} }{\phi^{2q}_{2i_1+2}(X)|X|}dX\right|^{\f{1}{2q}}\\
\les&\de^{\nu}K_1^{2q-1}K_0e^{-2q(\f{1}{2}-\ka)s}.
\end{align*}

For $I_2$, we use integration by parts and get
\begin{align*}
\int_{\mathbb{R}}\f{ w^{2q-1} }{\phi^{2q}_{2i_1+2}(X)|X|}\ve\p_XwdX=-\f{1}{2q}\int_{\mathbb{R}}\f{ w^{2q} }{\phi^{2q}_{2i_1+2}(X)|X|}\p_X\ve dX-\f{1}{2q}\int_{\mathbb{R}} w^{2q}\ve \p_X\left(\f{1}{\phi^{2q}_{2i_1+2}(X)|X|}\right)dX.
\end{align*}
Then by (\ref{vebd}), (\ref{pvebd}) and (\ref{phi}), it follows that
$$\left|\int_{\mathbb{R}}\f{ w^{2q-1}}{\phi^{2q}_{2i_1+2}(X)|X|}\ve \p_X w dX \right|\les e^{-(\f{1}{2}-\ka)s}\int_{\mathbb{R}}\f{ w^{2q} }{\phi^{2q}_{2i_1+2}(X)|X|}dX.$$

$I_3$ is similar to $I_2$,
$$\left|\int_{\mathbb{R}}\f{ w^{2q-1} }{\phi^{2q}_{2i_1+2}(X)|X|}w\p_X\ve dX \right|\les e^{-(\f{1}{2}-\ka)s}\int_{\mathbb{R}}\f{ w^{2q} }{\phi^{2q}_{2i_1+2}(X)|X|}dX.$$

For $I_4$, we first compute 
\begin{align*}
[A,\p_X]\ve=&-\left(\al_1+\p_X \sum_{l=1}^{L}e^{(\al_1-\al_l)s}\Psi_{i_l}(\bar{X}_l-\bar{Y}_l) \right)\p_X\ve=-\f{\al_1+\p_X \sum_{l=1}^{L}e^{(\al_1-\al_l)s}\Psi_{i_l}(\bar{X}_l-\bar{Y}_l)}{\al_1X+ \sum_{l=1}^{L}e^{(\al_1-\al_l)s}\Psi_{i_l}(\bar{X}_l-\bar{Y}_l)}A\ve,
\end{align*}
so by (\ref{X-Y}), (\ref{psil}) and (\ref{pxpsi}), we have
$$|[A,\p_X]\ve| \les\left|\f{A\ve}{X(\al_1+\f{ \sum_{l=1}^{L}e^{(\al_1-\al_l)s}\Psi_{i_l}(\bar{X}_l-\bar{Y}_l)}{X})}\right|\les \f{|A\ve|}{|X|}.$$
In this way, 
\begin{align*}
|I_4|\les\left|\int_{\mathbb{R}}\f{ w^{2q-1} }{\phi^{2q}_{2i_1+2}(X)|X|}\f{A\ve}{X}\ve\right| dX\les \|\f{\ve}{X}\|_{L^\infty}\int_{\mathbb{R}}\f{ w^{2q} }{\phi^{2q}_{2i_1+2}(X)|X|}dX\les e^{-(\f{1}{2}-\ka)s}\int_{\mathbb{R}}\f{ w^{2q} }{\phi^{2q}_{2i_1+2}(X)|X|}dX.
\end{align*}

For $I_5$, direct computation shows
\begin{align*}
&A\left(e^{(\al_1-\al_{l_1})s}\Psi_{i_{l_1}}(\bar{X}_{l_1}-\bar{Y}_{l_1})\p_{\bar{X}_{l_2}}\Psi_{i_{l_2}}(\bar{X}_{l_2}-\bar{Y}_{l_2})\right)\\
=&\left(\al_1X+\sum_{l=1}^{L}e^{(\al_1-\al_l)s}\Psi_{i_l}(\bar{X}_l-\bar{Y}_l)\right)\p_{\bar{X}_{l_1}}\Psi_{i_{l_1}}(\bar{X}_{l_1}-\bar{Y}_{l_1})\p_{\bar{X}_{l_2}}\Psi_{i_{l_2}}(\bar{X}_{l_2}-\bar{Y}_{l_2})\\
&+e^{(\al_1-\al_{l_1})s}\Psi_{i_{l_1}}(\bar{X}_{l_1}-\bar{Y}_{l_1})A\p_{\bar{X}_{l_2}}\Psi_{i_{l_2}}(\bar{X}_{l_2}-\bar{Y}_{l_2}),
\end{align*}
by (\ref{pxpsi}), we have 
\begin{align*}
&\left|\left(\al_1X+\sum_{l=1}^{L}e^{(\al_1-\al_l)s}\Psi_{i_l}(\bar{X}_l-\bar{Y}_l)\right)\p_{\bar{X}_{l_1}}\Psi_{i_{l_1}}(\bar{X}_{l_1}-\bar{Y}_{l_1})\p_{\bar{X}_{l_2}}\Psi_{i_{l_2}}(\bar{X}_{l_2}-\bar{Y}_{l_2})\right|\\
\les&\de^{\f{1}{2i_1}}e^{\al_1s}e^{(\f{\al_1}{\al_{l_1}}-1-\f{1}{\al_{l_1}})s}  e^{(\f{\al_1}{\al_{l_2}}-1-\f{1}{\al_{l_2}})s}\\
\les&\de^{\f{1}{2i_1}}e^{(\al_1+\f{\al_1}{\al_{l_1}}+\f{\al_{1}}{\al_{l_2}}-2-\f{1}{\al_{l_1}}-\f{1}{\al_{l_2}})s}.
\end{align*}
And by (\ref{psil}), (\ref{Apsi}),
\begin{align*}
&e^{(\al_1-\al_{l_1})s}\Psi_{i_{l_1}}(\bar{X}_{l_1}-\bar{Y}_{l_1})A\p_{\bar{X}_{l_2}}\Psi_{i_{l_2}}(\bar{X}_{l_2}-\bar{Y}_{l_2})\\
\les&\de^{\f{1}{2i_m(2i_m+1)}}e^{(\al_1-1)s}\de^{\f{1}{2i_1}}e^{(\al_{l_2}-\al_1+1)(\f{1}{2i_{l_2}+1}-2)s+\al_{l_2}s}\\
\les&\de^{\f{1}{2i_m(2i_m+1)}+\f{1}{2i_1}}e^{(\al_{l_2}-\al_1+1)(\f{1}{2i_{l_2}+1}-2)s+\al_{l_2}s+(\al_1-1)s}.
\end{align*}
Then Young's inequality implies that 
\begin{align*}
|I_5|\leq&\f{c_0}{q}\int_{\mathbb{R}}\f{ w^{2q} }{\phi^{2q}_{2i_1+2}(X)|X|}dX+C(q)\int_{\mathbb{R}}\f{(\de^{\f{1}{2i_m(2i_m+1)}+\f{1}{2i_1}}e^{\td{\al}s})^{2q}}{\phi^{2q}_{2i_1+2}(X)|X|}dX   \\
\leq&\f{c_0}{q}\int_{\mathbb{R}}\f{ w^{2q} }{\phi^{2q}_{2i_1+2}(X)|X|}dX+C(q)\de^{\f{q}{i_m(2i_m+1)}+\f{1}{i_1}}e^{(2q\td{\al}+\al_1)s}.
\end{align*}

For $I_6$, by (\ref{hl}), (\ref{Apsi})
\begin{align*}
&A\sum_{l=1}^{L}(\f{d}{ds}Y_l-\al_1 Y_l)\p_{\bar{X}_l}\Psi_{i_l}(\bar{X}_l-\bar{Y}_l)\\
\les&e^{\f{1}{2i_1}s}\de^{\f{1}{2i_m(2i_m+1)}}\de^{\f{1}{2i_1}}e^{(\al_l-\al_1+1)(\f{1}{2i_l+1}-2)s+\al_ls}\\
\les&\de^{\f{1}{2i_m(2i_m+1)}+\f{1}{2i_1}}e^{(\al_l-\al_1+1)(\f{1}{2i_l+1}-2)s+\al_ls+\f{1}{2i_1}}\\
\les&\de^\nu e^{\td{\al}s}.
\end{align*}
Similarly, 
$$|I_6|\leq\f{c_0}{q}\int_{\mathbb{R}}\f{ w^{2q} }{\phi^{2q}_{2i_1+2}(X)|X|}dX+C(q)\de^{\f{q}{i_m(2i_m+1)}+\f{1}{i_1}}e^{(2q\td{\al}+\al_1)s}.$$

At last, 
\begin{align*}
&\f{d}{ds}\left(\f{1}{2q} \int_{\mathbb{R}}\f{ w^{2q} }{\phi^{2q}_{2i_1+2}(X)|X|}dX  \right)\\
\leq&\left( \f{1}{2}-\f{C}{q}  \right)\int_{\mathbb{R}}\f{ w^{2q} }{\phi^{2q}_{2i_1+2}(X)|X|}dX+(\f{2c_0}{q}+3e^{-(\f{1}{2}-\ka)s})\int_{\mathbb{R}}\f{ \ve^{2q} }{\phi^{2q}_{2i_1+2}(X)|X|}dX\\
&+\de^{\td{p}}K_1^{2q-1}K_0e^{-2q(\f{1}{2}-\ka)s}+C\de^{\f{q}{i_m(2i_m+1)}+\f{1}{i_1}}e^{(2q\td{\al}+\al_1)s}\\
\leq&\left( \f{1}{2}-\f{C}{q}  \right)\int_{\mathbb{R}}\f{ \ve^{2q} }{\phi^{2q}_{2i_1+2}(X)|X|}dX+\de^{\nu}K_1^{2q-1}K_0e^{-2q(\f{1}{2}-\ka)s}+C\de^{\nu}e^{(2q\td{\al}+\al_1)s}.
\end{align*}
Take $|C/q|\leq\ka$, $p_1=2q\nu+\f{1}{2i_1}$, $p_2=2q\td{\al}+\al_1$ and integrate it,
\begin{align*}
&\int_{\mathbb{R}}\f{ w^{2q} }{\phi^{2q}_{2i_1+2}(X)|X|}dX\\
\leq&
e^{-2q(\f{1}{2}-\ka)(s-s_0)}\int_{\mathbb{R}}\f{ w_0^{2q} }{\phi^{2q}_{2i_1+2}(X)|X|}dX+Ce^{-2q(\f{1}{2}-\ka)s}\de^{p_1-p_2+(\f{1}{2}-\ka)}+\de^{\nu}K_1^{2q-1}K_0e^{-2q(\f{1}{2}-\ka)s}\\
\leq&e^{-2q(\f{1}{2}-\ka)s}(1+C\de^{p_1-p_2+(\f{1}{2}-\ka)}+\de^{\nu}K_1^{2q-1}K_0)\\
\leq&e^{-2q(\f{1}{2}-\ka)s}\f{K_1^{2q}}{2^{2q}},
\end{align*}
as long as $K_1\geq 2^{2q}(\f{1}{K_1^{2q-1}}+\de^{\nu}K_0+\f{C\de^{p_1-p_2+(\f{1}{2}-\ka)}}{K_1^{2q-1}}).$

\qed

\proof[Proof of Theorem 4.1]
Assume the conclusion is not true, which means there exists $\bar{\ve}$ whose initial data supported in $[\f{4}{3\de}, \f{5}{3\de}]$ and satisfying (\ref{ininorm}) is not global, whose lifespan of $\bar{\ve}$ is $s_*$. But \textbf{Proposition 4.3} implies that if $\bar{\ve}(s)$ is starting from $s_*$, $\|\bar{\ve}(s)\|_{\mathcal{B}}$ will keep small for a short time interval $[s_*,s_*+\eta),$ so does $\p_X\ve$ by \textbf{Corollary 4.4}. Then according to the continuation criterion of hyperbolic equation, $\bar{\ve}$ would not blow up in this time interval, which contradicts with former claim.

\qed

\end{document}